\documentclass{article}
\usepackage{graphicx,color}
\usepackage{amssymb}
\usepackage[geometry]{ifsym}
\usepackage{rotating}

\newtheorem{lemma}{Lemma}
\newtheorem{theorem}{Theorem}

\newcommand{\DONOTTEX}[1]{}

\begin{document}

\title{{\bf A note on uniquely \boldmath$10$-colorable graphs}}
\author{Matthias Kriesell}

\maketitle

\begin{abstract}
  \setlength{\parindent}{0em}
  \setlength{\parskip}{1.5ex}

  We prove for $k \leq 10$, that every graph of chromatic number $k$ with a unique $k$-coloring admits a clique minor of order $k$.


  
\end{abstract}

\maketitle

A {\em clique minor} of a (simple, finite, undirected) graph $G$ is a set of connected, nonempty, pairwise disjoint, pairwise adjacent subsets of $V(G)$,
where a set $A \subseteq V(G)$ is {\em connected} if $G[A]$ is connected,
and disjoint $A,B \subseteq V(G)$ are {\em adjacent} if there exists an edge $xy \in E(G)$ with $x \in A$ and $y \in B$.
An {\em anticlique} of $G$ is a set of pairwise nonadjacent vertices, and 
a {\em Kempe-coloring} of a graph $G$ is a partition $\mathfrak{C}$ into anticliques such that any two of them induce a connected subgraph in $G$.
In particular, if $A,B$ are distinct members of $\mathfrak{C}$ then every vertex from $A$ must have a neighbor in $B$ (*).
We recall the following facts from \cite{Kriesell2016}.

\begin{lemma} \cite{Kriesell2016}
\label{L1}
Every graph $G$ with a Kempe-coloring of order $k$ satisfies
$|E(G)| \geq (k-1)|V(G)|-{k \choose 2}$, with equality if and only if every pair of members of every Kempe-coloring of order $k$ induces a tree.
\end{lemma}

{\bf Proof.}
Let $\mathfrak{C}$ be a Kempe-coloring of order $k$ of $G$ and $A \not= B$ from $\mathfrak{C}$;
then $|E(G[A \cup B])| \geq |A|+|B|-1$ since $G[A \cup B]$ is a connected graph on $|A|+|B|$ vertices,
with equality if and only if $G[A \cup B]$ is a tree. Since $G[A \cup B]$ and $G[A' \cup B']$ are edge-disjoint for $\{A,B\} \not= \{A',B'\}$ we
get $|E(G)|=\sum |E(G[A \cup B])| \geq \sum (|A|+|B|-1)$, where the sums are taken over all subsets $\{A,B\}$ of $\mathfrak{C}$ with $A \not= B$.
Since every $X \in \mathfrak{C}$ occurs in exactly $k-1$ of these sets, the latter sum equals $(k-1)|V(G)|-{k \choose 2}$,
with equality if and only if any two members of $\mathfrak{C}$ induce a tree, which proves the statement for $\mathfrak{C}$.
As the latter bound is independent from the actual $\mathfrak{C}$,
equality holds for $\mathfrak{C}$ if and only if it holds for {\em all} Kempe-colorings of order $k$, which proves the Lemma.
\hspace*{\fill}$\Box$

\begin{lemma} \cite{Kriesell2016}
\label{L2}
Every graph with a Kempe-coloring of order $k$ is $(k-1)$-connected.
\end{lemma}

{\bf Proof.}
Let $\mathfrak{C}$ be a Kempe-coloring of order $k$ of a graph $G$. Then $|V(G)|>k-1$. 
Suppose, to the contrary, that there exists a separating vertex set $T$ with $|T|<k-1$. Then there exist $A \not= B$ in $\mathfrak{C}$ with
$(A \cup B) \cap T= \emptyset$; since $G[A \cup B]$ is connected, $A \cup B \subseteq V(C)$ for some component $C$ of $G-T$.
Now take any $x \in V(G) \setminus (T \cup V(C))$. Then $x$ is contained in some $Z \in \mathfrak{C}$ distinct from $A$ (and $B$),
but, obviously, $x$ cannot have a neighbor in $A$, contradicting (*).
\hspace*{\fill}$\Box$

An {\em $(H,k)$-cockade} is any graph isomorphic to $H$ and any graph that can be obtained by taking the union of two $(H,k)$-cockades
whose intersection is a complete graph on $k$ vertices. 
The following is the main result from \cite{SongThomas2006}.

\begin{theorem} \cite{SongThomas2006} 
\label{T0}
Every graph with $n>8$ vertices and at least $7n-27$ edges has a clique minor of order $9$, unless
it is isomorphic to $K_{2,2,2,3,3}$ or a $(K_{1,2,2,2,2,2},6)$-cockade.
\end{theorem}

Now we are prepared to prove the main statement of this note.

\begin{theorem}
\label{T1}
Every graph with a Kempe-coloring of order $10$ has a clique minor of order $10$.
\end{theorem}

{\bf Proof.}
Let $A \not= B$ be two color classes of a Kempe-coloring $\mathfrak{C}$ of order $10$ of a graph $G$.
Then $\mathfrak{C'}:=\mathfrak{C} \setminus \{A,B\}$ is a Kempe-coloring of $G':=G-(A \cup B)$, of order $8$.
By Lemma \ref{L1}, $G'$ is a graph on $n' \geq 8$ vertices with at least $7n'-28$ edges.

If $n'=8$ then $V(G')$ is a clique of order $8$, and, for every $x \in V(G')$, $G[\{x\} \cup A]$
and $G[\{x\} \cup B]$ are stars centered at $x$; therefore, if $ab$ is any edge in $G[A \cup B]$,
$V(G') \cup \{a,b\}$ is a clique of order $10$. So we may assume that $n' \geq 9$.

Now let $z$ be an endvertex of any spanning tree of $G[A \cup B]$.
Without loss of generality, we may assume that $z \in A$, otherwise we swap the roles of $A,B$.
Every $C \in \mathfrak{C}'$ contains a neighbor $x_C$ of $z$ in $G$ by (*). If these eight vertices form a clique then one checks readily
that $\{\{x_C\}:\,C \in \mathfrak{C}'\} \cup \{\{z\},(A \cup B) \setminus \{z\}\}$ is a clique minor in $G$ of order $10$
(every vertex $x_C$ has a neighbor in $B \subseteq (A \cup B) \setminus \{z\}$ by (*)).
Therefore, we may assume that $z$ has two distinct nonadjacent neighbors $x,y$ in $V(G')$.

If $G'+xy$ has a clique minor $\mathfrak{K}$ of order $9$ then we may assume without loss
of generality that $x$ is contained in some member $Q$ of $\mathfrak{K}$,
as $G'+xy$ is connected. Consequently,
$(\mathfrak{K} \setminus \{Q\}) \cup \{Q \cup \{z\},(A \cup B) \setminus \{z\}\}$
is a clique minor of $G$ of order ten (no matter whether $Q$ contains $y$ or not).

Hence we may assume that $G'+xy$ has no clique minor of order $9$.
As $G'+xy$ has at least $n' \geq 9$ vertices and at least $7n'-27$ edges,
we know that $G'+xy$ is one of the exceptional graphs in Theorem \ref{T0}.
By Lemma \ref{L2}, $G'$ is $7$-connected.
Therefore, $G'+xy$ is $7$-connected;
consequently, it cannot be the union of two graphs on more than $6$ vertices each, meeting in less than $7$ vertices.
It follows that $G'+xy$ is isomorphic to either $K_{2,2,2,3,3}$ or $K_{1,2,2,2,2,2}$, and $n'=11$ or $n'=12$.
Let $\mathfrak{B}$ be the set of single-vertex-sets in $\mathfrak{C}'$.
From $n' \geq |\mathfrak{B}|+2(8-|\mathfrak{B}|)$
we infer $|\mathfrak{B}| \geq 16-n'$, and, as $G[P \cup Q]$ is a star centered at the only vertex from $P$
for all $P \in \mathfrak{B}$ and $Q \in \mathfrak{C}' \setminus \{P\}$, every vertex from $\bigcup \mathfrak{B}$ is adjacent to all others of $G'$.
Consequently, $G'$ --- and hence $G'+xy$ --- has at least $16-n' \geq 4$ many vertices adjacent to all others.
However, $K_{2,2,2,3,3}$ has no vertex adjacent to all others, and $K_{1,2,2,2,2,2}$ has only one, a contradiction, proving the Theorem.
\hspace*{\fill}$\Box$

We may replace $10$ in Theorem \ref{T1} by any nonnegative $k<10$: Suppose that $G$ has a Kempe-coloring $\mathfrak{C}$ of order $k$ and
consider the graph $G^+$ obtained from $G$ by adding new vertices $a_{k+1},\dots,a_{10}$ and all edges from $a_i$, $i \in \{k+1,\dots,10\}$
to any other vertex $x \in V(G) \cup \{a_{k+1},\dots,a_{10}\}$.
Then $\mathfrak{C}^+:=\mathfrak{C} \cup \{\{a_{k+1}\},\dots,\{a_{10}\}\}$ is a Kempe-coloring
of $G^+$ of order $10$. By Theorem 1, $G^+$ has a clique minor $\mathfrak{K}$, and, as every $a_i$
is contained in at most one member of $\mathfrak{K}$,
the sets of $\mathfrak{K}$ not containing any of $a_{k+1},\dots,a_{10}$ form a clique minor of order at least $k$ of $G$.

A {\em $k$-coloring} of $G$ is a partition of $V(G)$ into at most $k$ anticliques, and the {\em chromatic number} $\chi(G)$ is the minimum number $k$
so that $G$ admits a $k$-coloring. (Observe that if a graph $G$ has a unique $k$-coloring then it has no $(k-1)$-coloring unless
it is a complete graph on less than $k$ vertices, so that, up to these exceptions, $\chi(G)=k$.)

\begin{theorem}
\label{T2}
For $k \leq 10$, every graph of chromatic number $k$ with a unique $k$-coloring admits a clique minor of order $k$.
\end{theorem}

{\bf Proof.}
Let $\mathfrak{C}$ be the unique $k$-coloring of $G$. Then $\mathfrak{C}$ is a Kempe-coloring of order $k$ (cf. \cite{Kriesell2016}),
and the statement follows from Theorem \ref{T1}.
\hspace*{\fill}$\Box$

{\bf Author's address:}

Department of Mathematics TU Ilmenau \\ Weimarer Stra{\ss}e 25, 98693 Ilmenau, Germany.

\end{document}